\documentclass[11pt]{article}
\newcommand{\ba}{\noindent $\begin{array}}
\newcommand{\ea}{\end{array}$}
\newcommand{\be}{\begin{equation}}
\newcommand{\ee}{\end{equation}}
\newcommand{\bd}{\begin{displaymath}}
\newcommand{\ed}{\end{displaymath}}
\newcommand{\beq}{\begin{eqnarray*}}
\newcommand{\eeq}{\end{eqnarray*}}
\newcommand{\beqn}{\begin{eqnarray}}
\newcommand{\eeqn}{\end{eqnarray}}

\usepackage{amsfonts,latexsym,lscape}

\small\normalsize

\usepackage{rawfonts}
\newfont{\Bb}{msbm10 scaled\magstep1}

\begin{document}

\pagestyle{plain}
\thispagestyle{empty}
\begin{center}
               { \textbf{Special Report for \\
         3rd INTERNATIONAL CONFERENCE ON THE ABS ALGORITHMS\\
                    May 13-14/01, Beijing}}
\end{center}
\vspace{2cm}
\begin{center}
{
              \large \textbf{VC++  ABSDLL V01 B}\\
\vspace{8pt}
      \LARGE \textbf{C++ CODES OF IMPLICIT LU ALGORITHMS FOR ABSDLL01 }
\vspace{1cm}

       \large          Xing Li, Ying Liu  and Antonino Del Popolo \\
\vspace{4cm}
 \vspace{7cm}
                   Centre for Optimization Research and Applications\\
        Department of Applied Mathematics, Dalian University of Technology
        }
\end{center}

\newpage
\setcounter{page}{1}

\begin{center}

{

}
\end{center}
\begin{center}
{\large \textbf{C++ CODES OF IMPLICIT LU ALGORITHMS FOR ABSDLL01
}}
\end{center}

\begin{center}
{\textbf{
                    Xing Li\footnote{CORA, Department of Applied Mathematics,
                                            Dalian University of Technology
                               Dalian 116024, China, xingli@student.dlut.edu.cn},
                   Ying Liu\footnote{Division of Software, Legend Computer
                               Company, Beijing 100085, China, liuyinge@legend.com.cn.}
                             } and \quad
                    Antonino Del Popolo\footnote{Department of
                                     Mathematics,
                                    University of Bergamo, Piazza Rosate 2,
                                    24129 Bergamo, Italy}
}
\end{center}

\vspace{2mm} \vspace{2mm}
\begin{center}
\parbox{13.5cm}{\small
\textbf{Abstract:}
  This report is devoted to some C++ codes implementing the implicit
  LU class algorithms for solving linear determined, and undetermined
  systems with $n$ variables and $m$ equations. A main program used in part of
  the numerical test is given in the last section.\\[10pt]

                           }
\end{center}
\setcounter{equation}{0}
\bigskip
\section{Introduction}
\setcounter{equation}{0}
ABS methods were introduced by Abaffy, Broyden and Spedicato (1982/84) \cite{AbBS 82}, \cite{AbBS 84} originally
for solving systems of linear equations. The basic ABS class was later generalized to the
so-called scaled ABS class and, subsequently, applied to linear-squares, nonlinear equations and
optimization problems, see for instance, Abaffy and Spedicatio \cite{AbSp 89}, Spedicato \cite{Sped 97}, \cite{Sped 99},
Spedicato, Xia Z. and Zhang L. \cite{SpXZ 00} and Zhang  L., Xia Z. and Feng E. \cite{ZhXF 99}.
In this paper we first review the general scheme for solving linear determined or undetermined systems.

\bigskip
Let us consider the general linear systems, where rank($A$)  is
arbitrary,
\begin{center}
    $Ax=b$,\quad $A \in I\!\!R^{m \times n}$
\end{center}
or
\begin{center}
    $a_i^Tx=b_i$,\quad $i=1,\cdots,m$
\end{center}
where
\begin{center}
$x \in I\!\!R^m$,\quad $b \in I\!\!R^n$,\quad $m \leq n$ \quad and \quad
   $
    A=\left[
              \begin{array}{c}
                  a_i^T \\
                  \cdots \\
                  a_m^T
              \end{array}
     \right]
  $
\end{center}
The steps of the unscaled ABS class of algorithms are defined as follows:

\bigskip
\noindent \textbf{Basic (Unscaled) ABS Class of Algorithms}: \cite{AbBS
84}, Algorithm\textbf{1} \cite{AbSp 89}:

\hspace{8.6cm}pp. 21-22\\
\textbf{\underline{Begin Algorithm}}
\begin{description}
\item[\textbf {  (A) }]  Initialization.\\[1pt]
                     Give  an arbitrary vector $ x_1 \in  I\!\!R^n $,  and an arbitrarily nonsingular matrix $ H_1 \in R^{n,n} $.\\
                     Set $ i = 1 $ and iflag=0.
\vspace{-2mm}
\item[\textbf{  (B)}] Computer two quantities.\\[1pt]
                      Compute
                       $$
                          \begin{array}{l}
                                   s_i = H_i a_i \\
                                   \tau_i=\tau^Te_i =a_i^T x_i-b^T e_i
                          \end{array}
                       $$
\vspace{-2mm}
\item[\textbf{  (C)}] Check the compatibility of the system of linear equations.\\[1pt]
                    If $ s_i \not= 0 $ then goto (D). \\
                    If $ s_i = 0 $ and $ \tau_i = 0 $ then set
                      $$
                                    \begin{array}{l}
                                         x_{i+1} = x_i \\
                                         H_{i+1} = H _i
                                  \end{array}
                     $$
                  and goto (F), the $ i $-th equation  is a linear combination of the previous equations. Otherwise  stop, the
                  system has no solution.
\vspace{-2mm}
\item[\textbf{ (D)}] Computer the search vector $p_i\in I\!\!R^n$ by
                   \begin{equation}\label{F2.1}
                        p_i=H_i^Tz_i
                   \end{equation}
                  where $z_{i}$, the  parameter of Broyden, is arbitrary
                  save that
                   \begin{equation}\label{F2.2}
                          z_i^TH_ia_i\neq 0
                    \end{equation}
\vspace{-2mm}
\item[\textbf{ (E)}] Update the approximation of  the solution $ x_i $ by
                   \begin{equation}\label{F2.3}
                        x_{i+1} = x_i - \alpha _i p_i
                    \end{equation}
                where the stepsize $ \alpha _i $ is computed by
                   \begin{equation}\label{F2.4}
                        \alpha _i = \tau_i /a_i^T p_i
                   \end{equation}
                  If $i=m$ stop; $x_{m+1}$ solves the system.
\vspace{-2mm}
\item[\textbf{ (F)}] Update the (Abaffian) matrix $ H_i. $
                   Compute
                   \begin{equation}\label{F2.5}
                        H_{i+1} = H_i - H_ia_iw_i^TH_i /w_i^TH_ia_i
                   \end{equation}
                   where $ w_i \in I\!\!R^n $, the parameter of Abaffy, is arbitrary save for the condition
                    \begin{equation}\label{F2.6}
                              w_i^T H_i a_i = 1 \textrm{ or } \neq 0
                     \end{equation}
\vspace{-2mm}
\item[\textbf{ (G)}] Increment the index $ i $ by one and goto (B).
\end{description}

\bigskip
\noindent
We define $n$ by $i$ matrices $A_i,\, W_i$
and $ P_i$ by
\begin{equation}\label{F2.7}
    \begin{array}{l}
                  A_i=(a_1,\,\cdots,\,a_i)^T, \quad
                  W_i=(w_1,\,\cdots,\,w_i),\quad
                  P_i=(p_1,\,\cdots,\,p_i)
   \end{array}
\end{equation}
\textbf{\underline{End Algorithm}}

\bigskip
\noindent
Some properties of the above recursion, see for instance, Abaffy and Spedicato (1989), \cite{AbSp 89},
are listed below that are the basic formulae for  use later on.
\begin{description}
\item[\textbf{ a. }]  Implicit factorization property
                     \begin{equation}\label{F1.8}
                     A_i^TP_i = L_i
                     \end{equation}
                    with  $ L_i $ nonsingular lower triangular.
\vspace{-8pt}
\item[\textbf{ b.}]  Null space characterizations
                   \begin{equation}\label{F1.9}
                   \begin{array}{c}
                   {\cal N} (H_{i+1}) ={\cal  R} (A_i^T),\quad
                   {\cal N} (H_{i+1}^T) ={\cal R} (W_i),\quad  \\[6pt]
                   {\cal N} (A_i) ={\cal R} (H_{i+1}^T)
                   \end{array}
                   \end{equation}
                 where ${\cal N}$= Null and ${\cal R}$=Range.
\vspace{-8pt}
\item[\textbf{ c.}]\ \ The linear variety containing all solutions to $Ax=b$
                   consists of the vectors of the form
                  \begin{equation}\label{F1.10}
                          x = x_{t+1} + H^T_{t+1}q
                 \end{equation}
                 where $ q \in I\!\! R^n $ is arbitrary.
\end{description}

Much progress on the computational aspect of ABS algorithms have been made since the ABS
algorithms were found in the early 1980's. In the last several years lots of work has been done on enlarging,
improving and completing ABSPACK that is a package of ABS algorithms in FORTRAN codes due
to Spedicato and his collaborators  \cite{Bodo 89}, \cite{Bodo 90a}, \cite{Bodo 00a}, \cite{Bodo 00b},
\cite{Bodo 01a},\cite{BoLS 00a}, \cite{BPLS 00a},  \cite{BPLS 00b},  \cite{BPLS 01}.

In the last two years we investigated possibility of establishing C++/VC++ ABS software, i.e.
ABS software written in C++/VC++. Part of work in the subject is presented in this paper concerning
the implicit LU algorithms.

In this report a text of some codes implementing the implicit LU algorithms for solving linear
determined and undetermined systems with $n$ variables and $m$ equations is given.
The codes are written in C++ language using Matrix class (with CMatrix as its class name)
and Vector class (with CVector as its class name) that are made by ourselves for constructing
the ABS software. The paper is organized  based on the following three algorithms:

\begin{enumerate}
\item  Function \emph{iLUa}: the implicit LU method, where $A$ is regular (i.e. all principal
                 submatrices are nonsingular), without pivoting for
                 determined and undetermined systems.
\item  Function \emph{iLUaPivotC}: the implicit LU method with column pivoting and explicit
             column interchanges for linear determined and undetermined systems. Since the column
             pivot cause that the ordering of components of solution is changed, it is necessary
             to recover the ordering of components of solution and related codes are given.
\item  Function \emph{iLUaPivotR}: the implicit LU method with row pivoting and explicit row
             interchanges for determined and undetermined systems.
\end{enumerate}

In the next section three schemes of implicit LU algorithms and some general properties are given.
Codes implementing these algorithms are listed in section 3. Finally a main program to test
algorithms is presented.

\bigskip
\section{Implicit LU Algorithms}
\bigskip
\subsection{Schemes of implicit LU algorithms}
We give a short description of three versions of the implicit LU algorithms below:

\textbf{Implicit LU  Algorithm Without Pivoting (iLUa)}: \cite{AbBS 82}, \cite{AbBS 84}, \cite{Sped 01}\\
(iLUa is abbreviation for implicit LU algorithm without pivoting, under the assumption \\
of $A$ is regular)\\
Set $ x_1=0 $, $ H_1=I $ and $ i=1 $\\
\textbf{For} $ i=1 $ to $m$ \textbf{do}
\vspace{-0.3cm}
\begin{description}
\item \hspace{3mm} Set \hspace{4mm} $ s_i = H_ia_i $
\vspace{-0.3cm}
\item \hspace{15mm} $ d_i = s_i^Te_i $
\vspace{-0.3cm}
\item \hspace{11mm} $ x_{i+1} = x_i-((a_i^Tx_i-b_i)/d_i)H_i^Te_i $
\vspace{-0.3cm}
\item \hspace{3mm} if $ i \le m $, then set
\vspace{-0.3cm}
\item \hspace{10mm} $ H_{i+1}=H_i-s_ie_i^TH_i/d_i $
\vspace{-0.3cm}
\end{description}
\textbf{enddo}

\bigskip
\noindent
\textbf{Implicit LU  Algorithm With Column Pivoting (iLUaPivotC)}: \cite{AbSp 89}, \cite{Bodo 00a}\\
(iLUaPivotC is the abbreviation of "implicit LU algorithm with column pivoting" and
column interchanges, without the regularity of $A$. The ordering of components of solution is changed because of the column pivot)\\
Set $ x_1=0 $, $ H_1=I $ and $ i=1 $\\
\textbf{For} $ i=1 $ to $m$ \textbf{do}
\vspace{-0.3cm}
\begin{description}
\item \hspace{3mm} Set $ s_i = H_ia_i $
\vspace{-0.3cm}
\item \hspace{3mm} (only $ (i-1)(n-i+1) $ nonzero elements of $ H_i$ are used)
\vspace{-0.3cm}
\item \hspace{3mm} Determine $ d_i=|s_i^Te_{k_i}| $, such that \quad $ |s_i^Te_{k_i}|=\max\{|s_i^Te_j|\,|\, j=i,\cdots,n\} $
\vspace{-0.3cm}
\item \hspace{3mm} (only $n-i+1$ nonzero elements of $s_i$ are used)
\vspace{-0.3cm}
\item \hspace{3mm} If $ k_i \neq i $, then swap columns of $A$ and elements of $x_i$ and $s_i$ with these indeices
\vspace{-0.3cm}
\item \hspace{3mm} Set \quad $ x_{i+1} = x_i-((a_i^Tx_i-b_i)/d_i)H_i^Te_i $
\vspace{-0.3cm}
\item \hspace{3mm} (only $i$ nonzero elements of $x_i$ are updated)
\vspace{-0.3cm}
\item \hspace{3mm} if $ i \le m $, then set \quad $ H_{i+1}=H_i-s_ie_i^TH_i/d_i $
\vspace{-0.3cm}
\item \hspace{3mm} (only $i(n-i)$ nonzero elements of $H_{i+1}$ are updated)
\end{description}
\vspace{-0.3cm}
\textbf{enddo}

\bigskip
\noindent
\textbf{Implicit LU Algorithm With Row Pivoting  (iLUaPivotR)}:\\
(iLUaPivotR stands for implicit LU algorithm with row pivoting and explicit row interchanges, without the
regularity of $A$ and change of components of solution vector)\\
Set $ x_1=0 $, $ H_1=I $ and $ i=1 $\\
\textbf{For} $ i=1 $ to $m$ \textbf{do}
\vspace{-0.3cm}
\begin{description}
\item \hspace{3mm} Determine $ d_i=|p_i^Ta_{k_i}| $,\, such that \quad $ |p_i^Ta_{k_i}|=\max\{|p_i^Ta_j|\,|\, j=i,\cdots,n\} $
\vspace{-0.3cm}
\item \hspace{3mm} If $ k_i \neq i $, then swap rows of $A$ and elements of $x_i$ and $b$ with these indices
\vspace{-0.3cm}
\item \hspace{3mm} Set $ s_i = H_ia_i $
\vspace{-0.3cm}
\item \hspace{3mm} (only $ (i-1)(n-i+1) $ nonzero elements of $ H_i$ are used)
\vspace{-0.3cm}
\item \hspace{3mm} Set \quad $ x_{i+1} = x_i-((a_i^Tx_i-b_i)/d_i)H_i^Te_i $
\vspace{-0.3cm}
\item \hspace{3mm} (only $i$ nonzero elements of $x_i$ are updated)
\vspace{-0.3cm}
\item \hspace{3mm} if $ i \le m $, then set \quad $ H_{i+1}=H_i-s_ie_i^TH_i/d_i $
\vspace{-0.3cm}
\item \hspace{3mm} (only $i(n-i)$ nonzero elements of $H_{i+1}$ are updated)
\end{description}
\vspace{-0.3cm}
\textbf{enddo}

\bigskip
\subsection{Properties of general implicit LU algorithms}
The general implicit LU algorithm is obtained by the parameter choices
$ H_1=I $, $ z_i=e_i $, $ w_i=e_i $. Some properties of this class of algorithms is listed as followings:
\begin{enumerate}
\item[(a)] The algorithm is well defined iff $A$ is regular (i.e., all principal submatrices
           are nonsingular). Otherwise pivoting has to be performed.
\item[(b)] Since $ W_i^TH_{i+1}=[I_i,0]^TH_{i+1}=0 $, the first $i$ rows
           of the Abaffian matrix $H_{i+1}$ must be zero. More precisely, the
           Abaffian matrix has the following structure, with $ K_i \in R^{n-i,i} $
            $$
                    H_{i+1}=\left[
                                    \begin{array}{cc}
                                    0     &   0 \\
                                    K_i   &   I_{n-i}
                                \end{array}
                        \right]
            $$
\item[(c)] Only the first $i$ components of $p_i$ can be nonzero and the $i$-th
           is unity. Hence the matrix $P_i$ is unit upper triangular, so that the
           implicit factorization $A=LP^{-1}$ is of the LU type, with unit on the diagonal.
\item[(d)] Only $K_i$ has to be updated. The algorithm requires $nm^2-2m^3$ multiplications
           plus lower-order terms. Hence, for $m=n$ there are $n^3/3$ multiplications plus
           low-order terms, which is the same cost as for the classical LU factorization or
           Gaussian elimination (which are two essentially equivalent process).
\end{enumerate}

\bigskip
\section{Codes of Implicit LU Algorithms}
\bigskip
\subsection{iLUa code}
\begin{verbatim}
function BOOL iLUa(CMatrix m_A, CVector v_b, CVector &v_x, double ep1,
                   double ep2)
\end{verbatim}
\tt // Under the condition of regularity of $A$, iLUa determines the solution\\
// of the linear systems $Ax=b$ ( $A$ with dimension $m \times n$, $m \leq n$), using \\
// implicit LU algorithm without pivoting\\
// the type of return value of function is BOOL, true if system has a \\
// solution, otherwise the return value is false\\
\begin{verbatim}
// m_A  = an object of CMatrix class, denote coefficient matrix
// v_b  = an object of CVector class, denote right hand-side vector
// v_x  = an object of CVector class, denote solution vector
//        (using reference operator to output the solution vector)
// ep1  = a double type value of dependency control parameter
// ep2  = a double type value of residual control parameter

{
   try
   {
      CVector v_x1;
      // declare an object of CVector class
      v_x1=v_x;

      CVector v_s;
\end{verbatim}
\vspace{-0.3cm}
\hspace{1cm} // declare an object of CVector class, denotes the vector $s_i=H_ia_i$
\vspace{-0.3cm}
\begin{verbatim}
      CMatrix m_H(1.0,m_A.m_cols);
      // declare an object of CMatrix class, denotes the Abaffian matrix,
      // the initial matrix is unit matrix
      double ns;
\end{verbatim}
\vspace{-0.7cm}
\begin{description}
\item \hspace{0.8cm} // declare a double precision number, denotes , $|| p ||$, the norm of $p$ \\
\vspace{-1cm}
\item \hspace{0.8cm} double r;\\
\vspace{-1cm}
\item \hspace{0.8cm} // declare a double precision number, $r_i=a_i^Tx_i-b_i$\\
\end{description}
\vspace{-1.5cm}
\begin{verbatim}
      int iflag=0;
      int i=1;

      // iteration
      while(i<=m_A.m_rows)
      {
\end{verbatim}
\vspace{-0.3cm}
\hspace{1.8cm} // compute $ s_i=H_ia_i $
\vspace{-0.3cm}
\begin{verbatim}
          if(i==1)
             v_s=m_A.GetRow(i).Trans();
          else
          {
             v_s[i-1]=0.0;
             for(int i1=i;i1<=m_A.m_cols;i1++)
             {
                double sum1=0.0;
                for(int j1=1;j1<i;j1++)
                   sum1=sum1+m_H.GetValue(i1,j1)*m_A.GetValue(i,j1);
                v_s[i1]=sum1+m_A.GetValue(i,i1);
             }
          }

          r=Dot(m_A.GetRow(i),v_x1)-v_b[i];
          ns=v_s.Module();
          if(ns<=ep1)
          {
             if(fabs(r)<=ep2)
\end{verbatim}
\vspace{-0.3cm}
\hspace{2.7cm}// $i$-th row is linearly dependent with first $i-1$ rows \\
\vspace{-1cm}
\begin{verbatim}
             {
                iflag=iflag+1;
                i=i+1;
             }
             else
             {
                // system is incompatible
                iflag=-i;
                AfxMessageBox("No Solution");
                break;
             }
          }
          else
          {
\end{verbatim}

\hspace{2.5cm} // update solution $x_i$\\
\vspace{-1cm}
\begin{verbatim}
             double temp=0;
             temp=r/v_s[i];
             for(int i2=1;i2<=i;i2++)
                v_x1[i2]=v_x1[i2]-temp*m_H.GetValue(i,i2);
\end{verbatim}
\vspace{-0.2cm}
\hspace{2.4cm} // update projection matrix $H_i$\\
\vspace{-0.8cm}
\begin{verbatim}
             if(i<m_A.m_cols)
             {
                for(int i3=i+1;i3<=m_A.m_cols;i3++)
                {
                   for(int j3=1;j3<=i;j3++)
                   {
                      double tt;
                      tt=(v_s[i3]*m_H.GetValue(i,j3))/v_s[i];
                      tt=m_H.GetValue(i3,j3)-tt;
                      m_H.SetValue(i3,j3,tt);
                   }
                }
                for(i3=1;i3<=i;i3++)
                   m_H.SetValue(i,i3,0.0);
             }
             i++;
          }
       }
       if(iflag>=0)
       {
          v_x=v_x1;
          return(true);
       }
       else
          return(false);
    }
    catch(CErrorException *e)
    {
       AfxMessageBox(e->GetErrorInfo());
       e->Delete();
       exit(1);
    }
}
\end{verbatim}
\subsection{iLUaPivotC code}
\begin{verbatim}
function BOOL iLUaPivotC(CMatrix m_A, CVector v_b, CVector &v_x, double ep1,
                        double ep2)
\end{verbatim}
\vspace{-0.3cm}
// iLUPivotC determines the solution of the linear systems $Ax=b$ ($A$ with\\
// dimension $ m \times n $, $ m \leq n $ ) using implicit LU algorithm with column pivoting\\
// and explicit column interchanges, and the order of components solution \\
// vector are changed\\
// the type of return value of function is BOOL, true if system has a solution\\
// otherwise the return value is false\\
\begin{verbatim}
// m_A  = an object of CMatrix class, denote coefficient matrix
// v_b  = an object of CVector class, denote right hand-side vector
// v_x  = an object of CVector class, denote solution vector
//        (using reference operator to output the solution vector)
// ep1  = double type value of dependency control parameter
// ep2  = double type value of residual control parameter

{
   try
   {
      CVector v_x1;
      // declare an object of CVector class
      v_x1=v_x;

      CVector v_s;
\end{verbatim}
\vspace{-0.3cm}
\hspace{1.2cm}// declare an object of CVector class, denotes the vector $s_i=H_ia_i$\\
\vspace{-0.9cm}
\begin{verbatim}
      CMatrix m_H(1.0,m_A.m_cols);
      // declare an object of CMatrix class, denotes the Abaffian matrix,
      // the initial matrix is unit matrix
      int ki;
      // declare an integer number, denotes the ordering of pivot
      double r;
      int iflag=0;
      double d;

      int* Index;
      // declare a pointer to type integer
      Index=new int[m_A.m_cols];
      // use new operator to allocate a pointer to an array of integer with
      // dimension m_A.m_cols
      for(int ii=1;ii<=m_A.m_cols;ii++)
         Index[ii]=ii;
      // initialize the array of integer
      int i=1;

      // iteration
      while(i<=m_A.m_rows)
      {
\end{verbatim}
\vspace{-0.3cm}
\hspace{1.7cm} // compute $s_i=H_ia_i$
\vspace{-0.2cm}
\begin{verbatim}
         if(i==1)
            v_s=m_A.GetRow(i).Trans();
         else
         {
            v_s[i-1]=0.0;
            for(int i1=i;i1<=m_A.m_cols;i1++)
            {
               double sum1=0.0;
               for(int j1=1;j1<i;j1++)
                  sum1=sum1+m_H.GetValue(i1,j1)*m_A.GetValue(i,j1);
                  v_s[i1]=sum1+m_A.GetValue(i,i1);
            }
         }
         // pivoting
         d=-1;
         for(int j2=i;j2<=m_A.m_cols;j2++)
         {
         if(d<fabs(v_s[j2]))
         {
            d=fabs(v_s[j2]);
            ki=j2;
         }
      }
      if(ki!=i)
\end{verbatim}
\vspace{-0.7cm}
\begin{description}
\item \hspace{0.8cm} // swap $i$-th column and $k_i$-th column of $A$
\vspace{-0.3cm}
\item \hspace{0.8cm} // swap $i$-th component and $k_i$-th component of $x_i$ and $s_i$
\end{description}
\vspace{-0.8cm}
\begin{verbatim}
      // save ordering of pivot in Index
      {
         m_A=m_A.SwapCol(ki,i);
         v_x=v_x.SwapElement(ki,i);
         v_s=v_s.SwapElement(ki,i);
         int it;
         it=Index[i];
         Index[i]=Index[ki];
         Index[ki]=it;
      }

      r=Dot(m_A.GetRow(i),v_x1)-v_b[i];
      if(fabs(v_s[i])<=ep1)
      {
         if(fabs(r)<=ep2)
\end{verbatim}
\vspace{-0.3cm}
\hspace{1.8cm}// $i$-th row is linearly dependent with first $i-1$ rows\\
\vspace{-1cm}
\begin{verbatim}
         {
            iflag=iflag+1;
            i=i+1;
         }
         else
           // system is incompatible
           {
              iflag=-i;
              AfxMessageBox("No Solution");
              break;
           }
       }
       else
       {
\end{verbatim}
\vspace{-0.3cm}
\hspace{2cm} // update solution $x_i$
\vspace{-0.3cm}
\begin{verbatim}
           double temp=0;
           temp=r/v_s[i];
           for(int i2=1;i2<=i;i2++)
              v_x1[i2]=v_x1[i2]-temp*m_H.GetValue(i,i2);
\end{verbatim}
\vspace{-0.2cm}
\hspace{2cm} // update projection matrix $H_i$ \\
\vspace{-1cm}
\begin{verbatim}
           if(i<m_A.m_cols)
           {
              for(int i3=i+1;i3<=m_A.m_cols;i3++)
              {
                 for(int j3=1;j3<=i;j3++)
                 {
                    double tt;
                    tt=(v_s[i3]*m_H.GetValue(i,j3))/v_s[i];
                    tt=m_H.GetValue(i3,j3)-tt;

                    m_H.SetValue(i3,j3,tt);
                 }
              }
              for(i3=1;i3<=i;i3++)
                m_H.SetValue(i,i3,0.0);
           }
           i++;
       }
       if(iflag>=0)
      // exchange components of solution to adapt the original system
      {
         for(int i4=1;i4<=m_A.m_cols;i4++)
         {
            for(int j4=1;j4<=m_A.m_cols;j4++)
            {
               if(Index[j4]==i4)
                  v_x[i4]=v_x1[j4];
            }
         }
         return(true);
      }
      else
         return(false);
   }
   catch(CErrorException *e)
   {
      AfxMessageBox(e->GetErrorInfo());
      e->Delete();
      exit(1);
   }
}
\end{verbatim}

\subsection{iLUaPivotR code}
\begin{verbatim}
function BOOL iLUaPivotR(CMatrix &m_A, CVector &v_b, CVector &v_x, double ep1,
                         double ep2)
\end{verbatim}
// iLUPivotC determines the solution of the linear systems $Ax=b$ ($A$ with \\
// dimension $ m \times n$, $ m \leq n$ ) using implicit LU algorithm with column pivoting\\
// and explicit column interchanges\\
// the type of return value of function is BOOL£¬true if system has a solution\\
// otherwise the return value is false\\
\begin{verbatim}
// m_A  = an object of CMatrix class, denote coefficient matrix
// v_b  = an object of CVector class, denote right hand-side vector
// v_x  = an object of CVector class, denote solution vector  (using reference
//        operator to output the solution vector)
// ep1  = double type value of dependency control parameter
// ep2  = double type value of residual control parameter

{
   try
   {
      CVector v_x1;
      // declare an object of CVector class
      v_x1=v_x;

      CVector v_s;
\end{verbatim}
\begin{description}
\vspace{-0.3cm}
\item \hspace{0.8cm} // declare an object of CVector class, denotes the vector $s_i=H_ia_i$
\vspace{-0.3cm}
\item \hspace{0.8cm} double r;
\vspace{-0.3cm}
\item \hspace{0.8cm} // declare a double precision number, $r_i=a_i^Tx_i-b_i$
\end{description}
\vspace{-0.3cm}
\begin{verbatim}
      CMatrix m_H(1.0,m_A.m_cols);
      // declare an object of CMatrix class, denotes the Abaffian matrix, the
      // initial matrix is unit matrix
      int ki;
      // declare an integer number, denotes the ordering of pivot
      double d;
      int iflag=0;
      int i=1;

      // iteration
      while(i<=m_A.m_rows)
      {
         // pivoting
         double mpt=0
         d=-1;
         for(int j2=i;j2<=m_A.m_rows;j2++)
         {
            mpt=Dot(m_H.GetRow(i),m_A.GetRow(j2));
            if(d<fabs(mpt))
         {
             d=fabs(mpt);
             ki=j2;
         }
      }
      if(ki!=i)
\end{verbatim}
\vspace{-0.7cm}
\begin{description}
\item \hspace{0.8cm} // swap $i$-th row and $k_i$-th row of $A$\\
\vspace{-1cm}
\item \hspace{0.8cm} // swap $i$-th component and $k_i$-th component of $x_i$ and $s_i$\\
\end{description}
\vspace{-1.5cm}
\begin{verbatim}
      // save ordering of pivot in Index
      {
         m_A=m_A.SwapRow(ki,i);
         v_x=v_x.SwapElement(ki,i);
         v_b=v_b.SwapElement(ki,i);
      }
      mpt=Dot(m_H.GetRow(i),m_A.GetRow(i));
\end{verbatim}
\vspace{-0.3cm}
\hspace{1cm} // compute $s_i=H_ia_i$
\vspace{-0.3cm}
\begin{verbatim}
      if(i==1)
         v_s=m_A.GetRow(i).Trans();
      else
      {
         v_s[i-1]=0.0;
         for(int i1=i;i1<=m_A.m_cols;i1++)
         {
            double sum1=0.0;
            for(int j1=1;j1<i;j1++)
               sum1=sum1+m_H.GetValue(i1,j1)*m_A.GetValue(i,j1);
            v_s[i1]=sum1+m_A.GetValue(i,i1);
         }
      }

      r=Dot(m_A.GetRow(i),v_x1)-v_b[i];
      if(fabs(mpt)<=ep1)
      {
         if(fabs(r)<=ep2)
\end{verbatim}
\vspace{-0.3cm}
\hspace{1.8cm}// $i$-th row is linearly dependent with first $i-1$ rows
\vspace{-0.3cm}
\begin{verbatim}
         {
            iflag=iflag+1;
            i=i+1;
         }
         else
         // system is incompatible
         {
            iflag=-i;
            AfxMessageBox("No Solution");
            break;
         }
      }
      else
      {
\end{verbatim}
\vspace{-0.3cm}
\hspace{1.8cm}// update solution $s_i$
\vspace{-0.3cm}
\begin{verbatim}
         double temp=0;
         temp=r/mpt;
         for(int i2=1;i2<=i;i2++)
            v_x1[i2]=v_x1[i2]-temp*m_H.GetValue(i,i2);
\end{verbatim}
\vspace{-0.3cm}
\hspace{2cm} // update projection matrix $H_i$
\vspace{-0.3cm}
\begin{verbatim}
           if(i<m_A.m_cols)
           {
              for(int i3=i+1;i3<=m_A.m_cols;i3++)
              {
                 for(int j3=1;j3<=i;j3++)
                 {
                    double tt;
                    tt=(v_s[i3]*m_H.GetValue(i,j3))/mpt;
                    tt=m_H.GetValue(i3,j3)-tt;

                    m_H.SetValue(i3,j3,tt);
                 }
              }
              for(i3=1;i3<=i;i3++)
                 m_H.SetValue(i,i3,0.0);
            }
            i++;
         }
      }
      if(iflag>=0)
      {
         v_x=v_x1;
         return(true);
      }
      else
         return(false);
   }
   catch(CErrorException *e)
   {
      AfxMessageBox(e->GetErrorInfo());
      e->Delete();
      exit(1);
   }
}
\end{verbatim}

\bigskip
\section{A Main Program for LU Algorithm}
\rm To make ABS algorithms be used more widely, we encapsulate Matrix class,\\
Vector class and the ABS algorithms modules into the library modules\\
ABSDLL.dll. A main program for the use of module(function) iLUaPivotC\\
that solves linear system $Ax=b$ (with Micchelli-Fiedler matrix as coefficient\\
matrix) is as following:

\begin{verbatim}
#include "stdafx.h"
#include "absalg.h"
#include "iostream.h"

void main()
{
    try
    {
        int m=1000;
        int n=1000;

        HINSTANCE hLib=AfxLoadLibrary("ABSDLL.dll");

        CMatrix a(m,n);
        CVector b(m);

        for(int i=1;i<=m;i++)
        {
           for(int j=1;j<=n;j++)
           {
               a.SetValue(i,j,abs(i-j));
           }
        }
        for(int k=1;k<=m;k++)
            b[k]=k;

        CVector xx(0.0,n);
        if(!iLUaPivotC(a,b,xx,1.0e-7,1.0e-7))
        {
            cout<<("Do you want a least-square solution?\endl");
            char c;
            cin>>c;
            if(c=='y'||c=='y')
               // calling the modules for solving the least-squares problem
            else
               exit(0);
        }
        else
        {
            CVector r;
            r=a*xx-b;
            er=r.Module();
            cout<<("errorbound=%f\n",er);

            AfxFreeLibrary(hLib);
        }
    }
    catch(CErrorException *e)
    {
        AfxMessageBox(e->GetErrorInfo());
        e->Delete();
    }
}
\end{verbatim}

\bigskip
\textbf{Remarks}\\
As for numerical experiments, we have made some tests, for example, matrices with elements $a_{ij}=|i-j|$,
 $1 \leq i\leq m$, $1 \leq j \leq n$ (Micchelli-Fiedler matrix) and matrices with elements $a_{ij}=|i-j|^2$,
 $1 \leq i\leq m$, $1 \leq j \leq n$. The results show that
the algorithms are efficient. Further numerical experiments are in progress.

\end{document}